\documentclass[12pt]{amsart}

\usepackage{xspace,amssymb,epsfig,euscript,xypic} 

\author{Ben Webster}

\title[Small linearly equivalent $G$-sets]{Small linearly equivalent
  $G$-sets and a construction of Beaulieu}

\curraddr{Department of Mathematics\\
Massachusetts Institute of Technology\\
77 Massachusetts Avenue\\
Cambridge, MA 02139}

\address{Department of Mathematics\\
  University of California, Berkeley\\ 
  Berkeley, CA 94720}

\email{bwebster@math.mit.edu} 

\thanks{This material is based upon work supported under a National
  Science Foundation Graduate Research Fellowship and partially
  supported by the RTG grant DMS-0354321. The LSU Research Experience
  Program is supported by the Louisiana Board of Regents Enhancement
  grant, LEQSF (1999-2001)-ENH-TR-17 and National Science Foundation
  grant, DMS-0097530}

\urladdr{http://math.mit.edu/\~{}bwebster/} 


\begin{document}

\newtheorem{thm}{Theorem}[section]
\newtheorem*{thm*}{Theorem}
\newtheorem{prop}[thm]{Proposition}
\newtheorem{lem}[thm]{Lemma}
\newtheorem{cor}[thm]{Corollary}
\newtheorem*{question}{Question}

\theoremstyle{remark}
\newtheorem{remark}{Remark}

\newcommand{\Aut}{\mathrm{Aut}}
\newcommand{\B}{\mathcal{B}}
\newcommand{\bfrac}[2]{\Big[\frac{#1}{#2}\Big]}
\newcommand{\bkh}{\backslash}
\newcommand{\bskip}{\baselineskip=2\baselineskip}
\newcommand{\C}{\mathbb{C}}
\newcommand{\Cbar}[2]{\overline{\C{#1}{#2}}}
\newcommand{\Cf}{\mathbb{C}}
\newcommand{\CF}[2]{\ensuremath{\mathfrak{C}(#1,#2)}}
\newcommand{\Cinf}{\EuScript{C}^{\infty}}
\newcommand{\cmp}{cyclic mod $p$\xspace}
\newcommand{\cp}{\fr c}
\renewcommand{\char}{{char}\,}
\newcommand{\CS}{\EuScript{CS}}
\newcommand{\deck}{\EuScript{D}}
\newcommand{\defl}[1]{\ind{}{#1}}
\newcommand{\eH}{[X_H]-[Y_H]}
\newcommand{\End}{\mathrm{End}}
\newcommand{\ES}[1]{\EuScript{#1}}
\newcommand{\Fix}{\mathrm{Fix}}
\newcommand{\fr}[1]{\mathfrak{#1}}
\newcommand{\Frat}{\mathrm{Frat}\,}
\newcommand{\Gal}[1]{\Gamma(#1 |\Q)}
\newcommand{\GrR}[1]{a(#1 G)}
\renewcommand{\H}{\mathcal{H}}
\newcommand{\id}{\mathrm{id}}
\newcommand{\ind}[2]{\mathrm{ind}^{#1}_{#2}}
\newcommand{\indp}[2]{\mathfrak{ind}^{#1}_{#2}}
\renewcommand{\inf}[1]{\res{}{#1}}
\newcommand{\inn}[1]{\langle #1\rangle}
\newcommand{\Iso}{\mathrm{Iso}}
\newcommand{\K}{\mathcal{K}}
\renewcommand{\ker}{\mathrm{ker}\,}
\renewcommand{\L}[1]{\mathfrak{L}(#1)}
\newcommand{\lap}[1]{\Delta_{#1}}
\newcommand{\lapM}{\Delta_M}
\newcommand{\lineq}{linearly equivalent\xspace}
\newcommand{\mc}[1]{\mathcal{#1}}
\newcommand{\mG}{m_G}
\newcommand{\mK}{m_{\K}}
\renewcommand{\O}{\mathcal{O}}
\newcommand{\Orb}{\mathrm{Orb}}
\newcommand{\pad}{\hat{\Z}_p}
\newcommand{\perm}[1]{\pi_{#1}}
\newcommand{\Q}{\mathbb{Q}}
\newcommand{\R}{\mathbb{R}}
\newcommand{\res}[2]{\mathrm{res}^{#1}_{#2}}
\newcommand{\resp}[2]{\mathfrak{res}^{#1}_{#2}}
\newcommand{\RG}{\EuScript{R}_G}
\newcommand{\rk}{\mathrm{rk}\,}
\newcommand{\V}[1]{\mathbf{#1}}
\newcommand{\vp}{\varphi}
\newcommand{\Stab}{\mathrm{Stab}}
\newcommand{\Sym}{\mathrm{Sym}}
\newcommand{\Tr}{\mathrm{Tr}}
\renewcommand{\V}[1]{\mathbf{#1}}
\newcommand{\Z}{\mathbb{Z}}
\newcommand{\Zp}{\Z/p}
\newcommand{\mindeg}[2]{\fr{md}_{#1}(#2)}
\newcommand{\CD}{\R[\Delta]}
\newcommand{\nc}{\newcommand}
 \nc{\set}[1]{\mathsf{set}_{#1}}
 \nc{\modu}[1]{#1\! -\! \mathbf{mod}}
 \nc{\Fp}{\mathbb{F}_p}
 \nc{\XM}{X_{\mc M}}
 \nc{\YM}{Y_{\mc M}}
 \nc{\lnq}{\overset{lin}{=}}
 \nc{\mdoo}{\frac{\mindeg \Q G}{\#G}}
\nc{\tG}{\tilde G}

\begin{abstract}
  Two $G$-sets ($G$ a finite group) are called linearly equivalent over
  a commutative ring $k$ if the permutation representations $k[X]$ and
  $k[Y]$ are isomorphic as modules over the group algebra $kG$.  Pairs
  of linearly equivalent non-isomorphic $G$-sets have applications in
  number theory and geometry. We characterize the groups $G$ for which
  such pairs exist for any field, and give a simple construction of
  these pairs. If $k$ is $\Q$, these are precisely the non-cyclic
  groups.  For any non-cyclic group, we prove that there exist $G$-sets
  which are non-isomorphic and \lineq over $\Q$, of cardinality $\leq
  3(\#G)/2$.
  
  Also, we investigate a construction of P. Beaulieu which allows us to
  construct pairs of transitive linearly equivalent $S_n$-sets from
  arbitrary $G$-sets for an arbitrary group $G$. We show that this
  construction works over all fields and use it construct, for each
  finite set $\mc P$ of primes, $S_n$-sets linearly equivalent over a
  field $k$ if and only if the characteristic of $k$ lies in $\mc P$.
\end{abstract}

\maketitle

\noindent
Let $G$ be a finite group, $X$ a $G$-set and $k$ a commutative ring.
The set of maps from $X$ to $k$, denoted $k[X]$, has a natural
structure of a $k$-module, and a natural action of $G$ given by
precomposition. That is, $k[X]$ is a module over the group algebra
$kG$. We call $k[X]$ the {\bf permutation representation of $G$ on
$X$}.

Two $G$-sets $X$ and $Y$ are called {\bf linearly equivalent} over $k$ if
$k[X]\cong k[Y]$ as $kG$-modules.  Linear equivalence is an
equivalence relation, and will be denoted $X\lnq_k Y$.
If no base ring is written, then it will be assumed to be $\Q$.

The definition of linear equivalence
was originally motivated by the following theorems from number theory:

\begin{thm*}[Perlis, \cite{rP77}]
Two number fields $E=\Q[\alpha]$ and $E'=\Q[\alpha']$ have identical
Dedekind zeta functions if and only if there is a Galois extension
$L/\Q$ containing $E, E'$ with Galois group $G=Gal(L/\Q)$, such that
$G\cdot\alpha\lnq_\C G\cdot\alpha'$.
\end{thm*}

\begin{thm*}[Perlis, Boltje, \cite{rP78,rB97}] Moreover, if
  $G\cdot\alpha\lnq_{\mathbb{F}_p} G\cdot\alpha'$, then the
  $p$-torsion subgroups of the class groups of $E$ and $E'$ are
  isomorphic.
\end{thm*}

Ideas along these lines have been further developed in
\cite{dS04,BB04}.  Similarly, pairs of linearly equivalent $G$-sets
have been used, first by Sunada in \cite{tS85}, in the construction of
pairs of manifolds with identical Laplacian and length spectra. By the
same yoga, they can be used to produce pairs of isospectral graphs
(see, for example, the work of Stark and Terras in \cite{hS00} on zeta
functions of graphs).

In order for two $G$-sets $X$ and $Y$ to be linearly equivalent, it is
necessary that $\#X=\#Y$, since $\#X=\dim k[X]$.  Thus, a new
invariant of the group $G$ and ring $k$ is the set
$$\deg_k(G)=\{n\ |\ \exists X,Y, \ X\lnq_k Y, \ X\ncong Y, \#X=\#Y=n\}$$

If $X\lnq_k Y$, then $X\sqcup\{*\}\lnq_k Y\sqcup\{*\}$ for any
singleton set $\{ * \} $ so
\begin{equation*}
\deg_k(G)=\{n\in \Z\ |\ n\geq \ell\}
\end{equation*}
for some integer $\ell$, which is the degree of the smallest pair of
\lineq and non-isomorphic $G$-sets. We call this integer $\mindeg
kG$.  If $\deg_k(G)$ is empty, we say $\mindeg k G=\infty$.  In this
paper we use a variety of group theoretic techniques to obtain
bounds on $\mindeg \Q G$.

In Section \ref{sec:basic-operations}, we recall the basic
operations of restriction and induction of $G$-sets, and explore
their interplay with linear equivalence. The relationship between
linear equivalence on different groups will a primary tool in this
paper.

In Section \ref{sec:an-existence-theorem}, we characterize those
groups for which there are pairs of non-isomorphic $G$-sets which
are linearly equivalent over a fixed field $k$.  In the case where
$k=\Q$, this is exactly the groups which are not cyclic.  While the
question of finding such $G$-sets which are transitive is quite
difficult and subtle (see \cite{bD00, wF80, rG83, rG85} for some
partial results), finding non-transitive examples is surprisingly
easy.  We then apply this construction to a number of special cases
where particularly simple pairs of \lineq $G$-sets can be found.

In Section \ref{sec:bounding-degrees}, we show, by an analysis of
cases, that
\begin{thm*}
  If $G$ is not cyclic, $\mindeg \Q G/\#G\leq 3/2$.
\end{thm*}\noindent
We give sharper bounds for some smaller classes of groups, including
the class of all non-solvable groups.

In Section \ref{sec:beaul-constr}, we discuss Beaulieu's
  construction \cite{B91} of pairs of transitive \lineq $S_n$-sets
starting from arbitrary pairs of \lineq $G$-sets.  We show that this
construction is independent of the field used, and obtain some
criteria for when the constructed sets are not isomorphic.

In Section
\ref{sec:an-application}, we apply this construction to find pairs of
$G$-sets linearly equivalent over any field whose characteristic lies
outside a given finite set of prime numbers.

\section*{Acknowledgments}
\label{sec:acknowledgments}

I would like to thank Robert Perlis, Bill Dunbar, Neal Stoltzfus and
Robert Snyder for their generous support and advice, and the Louisiana
State REU for a great opportunity to do mathematics as well as for their warm
hospitality.

\section{Basic Operations on $G$-sets}\label{sec:basic-operations}

The symbols $G,H$ and $K$ will denote finite groups throughout. By
convention, all $G$-actions are on the left.  We let $\set{G}$ be the
category of finite $G$-sets, with morphisms given by equivariant maps.
As usual, $\mathbb{F}_q$ denotes the finite field with $q$ elements,
considered as a field, or as an abelian group.

Throughout the rest of the paper, $k$ will always denote a field and
$p$ will always denote the characteristic of this field, which may 0
or a prime number.

If $X$ is a $G$-set, and $A$ any set, we define a new $G$-set $A\cdot
X$ to be the set $A\times X$ with $G$ acting on the right factor only
(we use $\cdot$ to distinguish this operation from the Cartesian
product of two $G$-sets, which has the diagonal action by definition).
We will denote $\{1,\ldots,n\}\cdot X$ by $n\cdot X$.

If $\psi:H\to G$ is a homomorphism, then there is a functor
$$\res{}\psi:\ \set{G}\to\set{H}$$
called {\bf restriction along} $\psi$, defined
 by composition of action
homomorphisms: if $X$ is a $G$-set with action homomorphism
$\rho_X:G\to\Sym X$, then $H$ acts by $\rho_X\circ \psi$. If the map
$\psi$ is injective, we think of $H$ as a subgroup of $G$, and
denote the restriction by $\res GH$.

We let $\inf\psi:\modu{kG}\to\modu{kH}$ denote the corresponding
functor for representations.

\begin{prop}
  For all $\vp:H\to G$,
  \begin{enumerate}
  \item The diagram
    \begin{equation*}
      \xymatrix@C=.5in{\set G\ar[r]^{\inf\psi{}}\ar[d]^{k[-]}& \set H\ar[d]^{k[-]}\\ \modu{kG}\ar[r]^{\inf\psi{}}&\modu{kH}}
    \end{equation*}
    commutes.
  \item \label{item:1} If $X\lnq_k Y$, then $\inf\psi X\lnq_k\inf\psi Y$.
  \item For $X$ any $G$-set, we have $\#X=\#\inf\psi X$.
  \end{enumerate}
\end{prop}

Part \eqref{item:1} of the above allows us to bound $\mindeg k G$ by
$\mindeg k {}$ of all its quotients.

\begin{lem}\label{L:quot}
If $\psi:H\to G$ is a surjective homomorphism, then $\mindeg k
H\leq\mindeg k G$.
\end{lem}
\begin{proof}
  Let $X$ and $Y$ be non-isomorphic $G$-sets \lineq over $k$ with
  $\#X=\mindeg kG$.  We know that $\inf\psi X\lnq_H \inf\psi Y$.
  Thus, we need only confirm that $\inf\psi X\ncong_H \inf\psi Y$.

  Since $X$ and $\inf\psi X$ have the same underlying set,
  any $H$-set isomorphism $\inf\psi X\to \inf\psi Y$ also defines an
  isomorphism $X\to Y$ (a priori this is only a set map).
  Obviously, this map commutes with the action of any group element in the
  image of $\psi$.  Since $\psi$ is surjective, this is in fact a
  $G$-set isomorphism.
  More elegantly, one could apply the functor $\defl\psi{}$ and use
  part~\eqref{item:8} of Proposition~\ref{induc-prop} below.

  Thus $\inf\psi X\ncong_H \inf\psi Y$, and $\mindeg kH\leq\#X=\mindeg kG$.
\end{proof}

Unfortunately, restriction is not useful for understanding $\mindeg
k G$ in terms of the subgroups of $G$.  Since most groups are richer
in subgroups than quotients, we will need an operation that can
transfer linearly equivalent $G$-sets to overgroups.

This is provided by the natural adjoint to the functor $\inf\psi$,
which is called {\bf induction along} $\psi.$ For any $H$-set $X$,
we let $\defl\psi X=(G\times X)/H$, where $H$ acts on the Cartesian
product $G\times X$ by
$$^h(g,x)=(g\psi(h^{-1}),{^hx}).$$
This defines a functor $\defl\psi{}:\set H\to\set G$.

We will also use $\defl\psi{}$ to denote the analogue of this
functor from representation theory: for a $kH$-module $V$, we define
$\defl\psi V=kG\otimes_{kH}V$, using the fact that $kG$ is a right
$kH$-algebra.

\begin{prop}\label{induc-prop}
  For all $\vp:H\to G$,
  \begin{enumerate}
  \item\label{item:2} The functor $\defl\psi{}$ is left adjoint to $\inf\psi{}$.
  \item\label{item:3}  The diagram
    \begin{equation*}
      \xymatrix@C=.7in{\set H\ar[r]^{\defl\psi{}}\ar[d]^{k[-]}& \set
      G\ar[d]^{k[-]}\\ \modu{kH}\ar[r]^{\defl\psi{}}&\modu{kG}}
    \end{equation*}
    commutes.
  \item \label{item:4}If $X\lnq_k Y$, then $\defl\psi X\lnq_k\defl\psi Y$.
  \item \label{item:5}For any subgroup $K\subset H$, $\defl\psi(H/K)\cong
    G/\psi(K)$. In particular,
  \begin{equation*}
    \#\defl\psi(H/K)=
    \frac{(G:\mathrm{im}\,\psi)}{(\ker\psi:K\cap\ker\psi)}\#(H/K).
  \end{equation*}
  \item  \label{item:6} If $\psi$ is injective, then $\#\ind GH X=(G:H)\#X$ for any $H$-set
    $X$.
  \item \label{item:7} If $\psi$ is injective, then the stabilizer of any element of
  $\defl\psi X$ is isomorphic to the stabilizer of some element of $X$.
  \item \label{item:8} If $\psi$ is surjective, then  $\defl\psi X\cong \ker\psi\backslash X$,
  the orbit space for the action of $\ker\psi$.  In particular, $\defl\psi(\inf\psi(X))\cong_G X$ for any $G$-set $X$.
  \end{enumerate}
\end{prop}
\begin{proof}
  Part \eqref{item:5}:  Under the action specified above $G\times H/K$ is an $G\times
  H$-set.  The stabilizer of the coset $(1,K)$ is $K'=\{(\psi(x),x)|x\in K\}$,
  and this action is transitive, so $G\times H/K\cong (G\times H)/K'$.
  Thus $\defl\psi H/K$ is transitive as a $G$-set, and the stabilizer
  of $H(1,K)$ is $\psi(K)$.
  
  The stabilizer of $H$ acting on the coset $(1,K)$ is $K\cap \ker\psi$.
  Furthermore, $G$-acts transitively on all $H$-orbits, so they are all
  of cardinality $(H:\ker\psi\cap K)$.  Thus, we calculate that
  \begin{align*}
      \#\defl\psi(H/K)&=\frac{\#
    G}{(H:\ker\psi\cap
    K)}\#(H/K)\\
&=\frac{\#G}{\#\mathrm{im}\,\psi(\ker\psi:\ker\psi\cap
    K)}\#(H/K).
  \end{align*}

  Part \eqref{item:6}:  Apply the formula of part~\eqref{item:5} to
  each component.

  Part \eqref{item:7}: This only needs to be checked for $X$ a
  transitive $G$-set.  In this case, $X\cong H/K$ for some subgroup
  $K\subset H$, and $\defl\psi X\cong G/\psi(K)$.  Thus the stabilizer
  of any element of $X$ is conjugate in $H$ to $K$, and thus
  isomorphic to $K$.  Similarly, the stabilizer of any element of
  $\defl\psi X$ is conjugate in $G$, and thus isomorphic to $\psi(K)$.
  Since $\psi$ is injective, the two stabilizers just be isomorphic.

  Part \eqref{item:8}: For each transitive $G$-set $G/K$,
  \begin{equation*}
      \defl\psi(\inf\psi(G/K))\cong_G\defl\psi(H/\psi^{-1}(K))\cong_GG/K.
  \end{equation*}
  Since $\inf\psi{}$ and $\defl\psi{}$ respect disjoint union, this
  implies the result for all $X$.
\end{proof}

\begin{cor}
  If $X\lnq_k Y$, then $X$ and $Y$ have the same number of orbits.
\end{cor}
\begin{proof}
  Let $\tau:G\to 1$ be the trivial homomorphism.  Then by part
  \eqref{item:8} of Proposition~\ref{induc-prop},  $\defl\tau
  X\cong G\backslash X$, the orbit space of $X$. By
  part \eqref{item:3} of the same proposition, $G\backslash X$ and
  $G\backslash Y$ are isomorphic as sets
  with an action of the trivial group, i.e. they are sets with the
  same cardinality.
\end{proof}
Unfortunately, one must be more careful when using induction than when
using restriction, since if we have non-isomorphic $H$-sets $X$ and $Y$,
it may still be that $\ind GHX\cong_G\ind GHY$ (unlike restriction along
a surjective map, induction is not full)

For instance, let $K_1,K_2\subset H$ are subgroups which are not
conjugate in $H$, but which are conjugate in $G$.  Of course,
$H/K_1\ncong_H H/K_2$, but using part \eqref{item:5} of
Proposition~\ref{induc-prop}, we see that
\begin{equation*}
  \ind GH(H/K_1)\cong_G G/K_1\cong_G G/K_2\cong_G \ind
GH(H/K_2),
\end{equation*}
since in $G$, the subgroups $K_1$ and $K_2$ are conjugate.

For example, if $A_4$ is the alternating group on 4 elements, and
$K_4$ the Klein four-group generated by the permutations
$\{(12)(34),(13)(24),(14)(23)\}$, then all elements of order two are
conjugate to each other in $A_4$, even though they are not in $K_4$.

However, if there is an element $x\in X$ such that for all $y\in Y$ the
stabilizers $\Stab_H(x)$ and $\Stab_H(y)$ are not isomorphic as abstract
groups, then $\defl\psi X\ncong_G \defl\psi Y$, since
induction along an injective map preserves the isomorphism class of
stabilizers, by Proposition~\ref{induc-prop}.

\begin{lem}\label{T:ind}
  If $\psi:H\to G$ is injective, if $X$ and $Y$ are non-isomorphic
  linearly equivalent $G$-sets, and if there is $x\in X$ such that for
  all $y\in Y$ the stabilizers $\Stab_H(x)$ and $\Stab_G(y)$ are not
  isomorphic, then $\mindeg kG\leq (G:H)\#X$.
\end{lem}

\section{An Existence Theorem}\label{sec:an-existence-theorem}

As is well known, the representations $\Q[X]\cong\Q[Y]$ are
isomorphic if and only if they have the same character. We denote
this character $\pi_X$. This character is know to be given by
$\pi_X(g)=\#\Fix_gX,$ the number of fixed points for the action of
$g$.

Since we will be interested in equivalence over all characteristics,
we will need an analogue of character over fields of positive
characteristic. In the case of permutation representations, there is
a very nice solution to this problem, in which needs only to
consider the fixed points of a more general class of groups than
cyclic ones.

We call a group \textbf{\cmp} if it is an extension of a cyclic
group by a $p$-group.  By convention, a $0$-group is trivial, so ``cyclic
mod 0'' simply means ``cyclic.''
\begin{lem}\label{L:restr}
  For two $G$-sets, $X$ and $Y\!$, the following are equivalent:
  \begin{enumerate}
  \item $X\lnq_k Y$.
  \item $\res GH X\cong_H\res GH Y$ for all \cmp subgroups $H\subset G$.
  \item $\#\Fix_X(H)=\#\Fix_Y(H)$ for all \cmp subgroups $H\subset G$.
\end{enumerate}
\end{lem}
In the case where $p=0$, this obviously reduces to the statement
that representations are determined by their characters.
\begin{proof}
$(1)\Rightarrow(2)$: Restriction is a functor, so it sends
isomorphisms to isomorphisms.

$(2)\Rightarrow(3)$:  Now, fix a \cmp subgroup $H\subset G$
and let $P$ be the unique Sylow $p$-subgroup and let
$F=\Fix_X(P)$. Thus $H/P$ is a cyclic $p'$-group.  Note that $F$ is naturally an
$H/P$-set, and the $H/P$ action on $k[F]$ is induced
from its embedding into $k[\res GH X]$.  

In this case, $k[F]$ is the
maximal trivial \emph{summand} of $\res GP k[X]$, because over
characteristic $p$, no nontrivial transitive permutation representation
of a $p$-group has any trivial summands.  This subspace is not unique,
since Krull-Schmidt decomposition is not unique, but its isomorphism
class as a $H/P$ representation is well-defined.  

We note that the classes of the permutation representations of $H/P$
are linearly independent in the representation ring of $H/P$ over $k$
(this result is true only for cyclic groups of order coprime to the
characteristic of the field).  Thus, $F$ can be reconstructed as a
$H/P$-space up to isomorphism from $k[F]$ as a $G/H$-set.

 Of course, 
\begin{equation*}
  \Fix_F(H/P)=\Fix_X(H),
\end{equation*}
so the number of fixed points of $H$ is determined by $k[\res GH X]$.

$(3)\Rightarrow(1)$: This is the hardest implication, and we
will not present a proof here.  See \cite[81.25 and 81.28]{cC87}.
\end{proof}

We will need a simple computation of $\#\Fix_{G/H}(K)$ for a subgroup
$K\subset G$.  First, we let
\begin{equation*}
L_K(H)=\{g\in G | K\subset gHg^{-1}\}.
\end{equation*}
Note that $L_K(H_1)\cap L_K(H_2)=L_K(H_1\cap H_2)$ and that
$L_K(H)H=L_K(H)$. That is, $L_K(H)$ is equipped with a free right
$H$-action by multiplication.

\begin{lem}\label{L:fixcount}
\begin{equation}\label{eq:perm1}
\#\Fix_{G/H}(K)=\#(L_K(H)/H)=\frac{\# L_K(H)}{\#H}
\end{equation}
\end{lem}
\begin{proof}
The natural map $L_K(H)/H\to G/H$ is an injection, and its image is $\Fix_{G/H}(K)$.
\end{proof}

Let $\mc M=\{M_1,\ldots,M_\ell\}$ be a set of
distinct subgroups of $G$ such that:
\renewcommand{\labelenumi}{(\arabic{enumi})}
\begin{enumerate}
\item \label{union} The union $H=\bigcup_{i}M_i$ is a subgroup of $G$.
\item \label{proper} For all $i$, $M_i$ is a proper subgroup of $H$.
\item \label{contain} Each \cmp subgroup $C\subset H$ is also contained in $M_i$ for some $i$.
\end{enumerate}

Note that in the case where $p=0$, the condition~\eqref{contain} is implied by condition~(\ref{union}).

If $G$ is not \cmp, then all the
proper maximal subgroups of $G$, or all the \cmp subgroups of $G$
will serve as such a set.

On the other hand, the subgroup $H$ must be non-\cmp, since otherwise,
condition~(\ref{contain}) implies $H=M_i$ for some $i$, which
condition~(\ref{proper}) explicitly forbids.  In particular, if $G$ is
\cmp, then no such subsets of the subgroups of $G$ exist.

We let $\mathcal{J}_e$ be the set of non-empty subsets of
$\{1,\ldots, \ell\}$ which have even cardinality, and
$\mathcal{J}_o$ be those of odd cardinality. Define $G$-sets $\XM$
and $\YM$
\begin{align}\label{eq:defXH}
  \XM&=\#H\cdot \frac GH \bigsqcup
    \left[\bigsqcup_{S\in\mathcal{J}_e}\#\left(\bigcap_{i\in S}M_i\right)\cdot
    \frac G{\bigcap_{i\in S} M_i}\right]\\ \label{eq:defYH}
  \YM&=\bigsqcup_{S\in\mathcal{J}_o}\#\left(\bigcap_{i\in S}M_i\right)\cdot
    \frac G{\bigcap_{i\in S} M_i}.
\end{align}
Note that $\XM\ncong\YM$, since $\Fix_{\XM}(H)$ is not empty, and, by condition~(\ref{proper})
$\Fix_{\YM}(H)$ is.

\begin{thm}\label{T:exist}
  Let $\mc M$ satisfy conditions
  (\ref{union}), (\ref{proper}) and (\ref{contain}) .
  Then for any field $k$ of characteristic $p$,
  \begin{equation*}
    \XM\lnq_k\YM.
  \end{equation*}
\end{thm}
\begin{proof}
Using formula \eqref{eq:perm1},
\begin{align}\label{eq:XHperm}
  \Fix_{\XM}(C)&=\#L_C(H)+\sum_{S\in\mathcal{J}_e}\#L_C\left(\bigcap_{i\in S}M_i\right)\\
  \Fix_{\YM}(C)&=\sum_{S\in\mathcal{J}_o}\#L_C\left(\bigcap_{i\in S}M_i\right)
\end{align}
By condition~(\ref{contain}), $C\leq M_i$ for some $i$.
By inclusion-exclusion,
\begin{equation}
\#L_C(H)=\sum_{S\in\mathcal{J}_o} \#\bigcap_{i\in
S}L_C\left(M_i\right) -\sum_{S\in\mathcal{J}_e} \#\bigcap_{i\in S}L_C\left(M_i\right).
\end{equation}
Substituting this into \eqref{eq:XHperm}, we find that
\begin{equation}
  \Fix_{\XM}(C)=\Fix_{\YM}(C)
\end{equation}
for all $g\in G$, so $\XM$ and $\YM$ are \lineq over any field of characteristic $p$.
\end{proof}

\begin{remark}
  If the set $\mc M$ contains all maximal subgroups of
  $H$, it leads to a new formula for the idempotents in the Burnside
  ring, which were described by Solomon in \cite{lS67}, and in this case our
  $G$-sets $X$ and $Y$ could also be defined using the formula given
  by Gluck in \cite{dG81}.

  On the other hand, our formula is more general, since it allows to
  choose smaller sets of subgroups, which will result in smaller
  $G$-sets.
\end{remark}

Combining Lemma~\ref{L:restr} and Theorem~\ref{T:exist}, we see that

\begin{cor}
There exist non-isomorphic, \lineq $G$-sets over $k$ if and only if $G$ is
not \cmp.
\end{cor}

\subsection{Frobenius Groups}
\label{sec:frobenius-groups}

In this subsection, we only consider the case where $p=0$ (for
example, $k=\Q$).

We call $G$ a {\em Frobenius group} if there is a
non-trivial\label{frobdef} proper subgroup $H\subset G$ such that
$H\cap gHg^{-1}=\{1\}$ for all $g\in G\backslash H$.  This is
equivalent to the existence of a transitive, non-regular $G$-set $X$
such that each non-trivial element fixes exactly 1 point, or none.

By Frobenius's theorem \cite[8.5.5]{dR96}, if $G$ is a Frobenius
group there is a normal subgroup $K\lhd G$ such that $K\cap gHg^{-1}
=\{1\}$ for all $g\in G$ and $G=KH$. In fact, $K$ is simply the
elements $G$ which are not conjugate to any nontrivial element of $H$.

Thus the set $\mc{F}=\{K, H, g_1Hg_1^{-1},\ldots\}$, where
$\{1,g_1,\ldots\}$ contains exactly one element from each coset of
$H$, satisfies $F_1\cap F_2=\{1\}$ for $F_1,F_2\in\mc{F}, F_1\neq
F_2$  and $\bigcup_{F\in\mc{F}}F=G$, i.e., $\mc{F}$ is a partition of
$G$ in the sense of \cite{rA70}.  Thus we may calculate
\begin{align*}
X_{\mc{F}}&=\#G\cdot\frac GG\bigsqcup(2^{\#\mc{F}-1}-1)\cdot\frac G{\{1\}}\\
Y_{\mc{F}}&=\left(\bigsqcup_{F\in\mc{F}}\#F\cdot \frac
GF\right)\bigsqcup(2^{\#\mc F-1}-
 (G:H)-1)\cdot \frac G{\{1\}}.
\end{align*}
Removing redundant copies of the regular action, we find the
$G$-sets
\begin{align}
\label{eq:Xfrob} \tilde{X}_{\mc{F}}&=\#G\cdot\frac
GG\bigsqcup(G:H)\cdot\frac G{\{1\}}\\ \label{eq:Yfrob}
\tilde{Y}_{\mc{F}}&=\bigsqcup_{F\in\mc{F}}\#F\cdot \frac
GF\cong_G\#K \cdot\frac GK\bigsqcup \#G\cdot\frac GH.
\end{align}
are \lineq.

Note that $\tilde{X}_{\mc F}\cong (G:H)\cdot X'_{\mc F}$ and
$\tilde{Y}_{\mc F}\cong (G:H)\cdot Y'_{\mc F}$, where
\begin{align*}
X'&= \# H\cdot\frac GG\bigsqcup \frac G{\{1\}}\\
Y'&= \# H\cdot\frac GH\bigsqcup \frac GK.
\end{align*}

This is, of course, considerably easier than the calculation we
would have to do for the $G$-sets $\XM$ and $\YM$ using all
maximal subgroups or all cyclic subgroups, and gives us much
smaller $G$-sets; $\#X'=\#Y'=\#G+\#H$.

Thus, if $G$ is a non-regular Frobenius group, $\mindeg \Q G\leq
\#G+\#H$.

In fact, we have
\begin{prop}
If $G$ is a non-regular Frobenius group, $$\frac{\mindeg \Q G}{\#
G}\leq \frac 43.$$
\end{prop}
\begin{proof}
From the discussion above, we know that $$\frac{\mindeg \Q G}{\#
G}\leq 1+\frac{1}{\#K}.$$  Since $khk^{-1}\in H$ if and only if $k$ or $h$
is the identity, the natural map $H\to\Aut K$ is injective.  Thus, $\#K\geq 3$.
\end{proof}

A computer search conducted with the computer algebra system GAP shows
that $\mindeg \Q{S_3}=8$, so this bound is strict.

Unfortunately, these $G$-sets will often not be equivalent over a
field of characteristic which divides the order of the group.

\subsection{$q$-groups}

In this subsection, $k$ is of characteristic $p$ which may be positive or
0.

We consider the case where $G$ is a $q$-group, for $q$
a prime different from $p$.  Since $q$-groups have so many quotients,
we can hope to find a quotient $\tilde G$ of any $q$-group $G$ which
is both complicated enough that $\mindeg k{\tilde G}$ is relatively
low, but its computation is tractable.  Luckily, such a quotient is
already provided by classical group theory.

First, we define the {\em Frattini subgroup} of $G$, $\Frat G$, to be
the intersection of all the maximal subgroups of $G$. This subgroup is
normal, since any conjugate of a maximal subgroup is maximal.  We
define $\tilde G=G/\Frat G$

Note that if $\{a_1,\ldots, a_n\}$, with $a_i\in G$ is a generating set
of $G$ if and only if its image in $\tilde{G}$ is as well.  Thus $G$ is
cyclic if and only if $\tilde{G}$ is.

For a general group, this quotient is rather hard to compute, but if $G$ is a
$q$-group, we can obtain important information about $\tilde G$ from
the Burnside Basis Theorem:
\begin{thm}\emph{(Burnside, \cite[5.3.2]{dR96})}\label{burn-base}
  If $G$ is a $q$-group, then $\tilde G=G/\Frat G$ is the largest
  elementary abelian quotient of $G$.  
\end{thm}

Thus using this reduction, we obtain the

\begin{prop}\label{P:pgroup} If $G$ is any non-cyclic $q$-group with
  $q\neq p$,
  \begin{equation*}
    \frac{\mindeg k G}{\# G}\leq \frac {q+1}q.
  \end{equation*}
  In particular, if $k=\Q$, this bound holds for all $q$.
\end{prop}
Of course, $\mindeg k G=\infty$ if $G$ is a $p$-group (since $G$ is
cyclic mod $p$).
\begin{proof}
  Since $G$ is non-cyclic, Theorem~\ref{burn-base} shows that $\tilde
  G\cong\mathbb{F}_q^n$, with $n\geq 2$.  Fix a subspace $K\subset
  \tG$ of codimension $2$ (that is $K\cong (\mathbb{F}_q)^{n-2}$). For
  any $g\in \tG$, $\langle g,K\rangle$ is a proper subspace, so every
  element of $\tG$ is in a proper subgroup containing $K$.  We let
  $\mc{A}$ be the set of maximal subgroups of $\tG$ containing $K$.
  
  This satisfies the first two hypotheses of Theorem~\ref{T:exist}.
  Since $p$ doesn't divide the order of $G$, any \cmp subgroup of $G$
  is actually cyclic.  Thus, the sets $X_{\mc A}$ and $Y_{\mc A}$
  described in Theorem~\ref{T:exist} are linearly equivalent over $k$.
  
  Now, if $A_1,A_2\in \mc{A}, A_1\neq A_2$, then $A_1+A_2=\tG$, by
  maximality, so $\dim(A_1\cap A_2)=n-2$, and since it contains $K$,
  $K=A_1\cap A_2$. Thus, removing isomorphic orbits from $X_{\mc A}$
  and $Y_{\mc A}$, and dividing out by $q^{n-2}$, we find that
\begin{align}
\label{eq:XF1}
X'_{\mc{A}}=\frac {\tG}K\bigsqcup q^2\cdot\frac {\tG}{\tG}\\
\label{eq:XF2}
Y'_{\mc{A}}=\bigsqcup_{A\in\mc{A}}q\cdot\frac{{\tG}}{A}
\end{align}
are \lineq over k.

Thus, $\frac{\mindeg k{\tG}}{\#\tG}\leq \frac{q+1}{q}$.  Applying Lemma \ref{L:quot}, we obtain the desired result.
\end{proof}

\section{Bounding Degrees}\label{sec:bounding-degrees}

In this section, we will only consider the case where $k$ is
characteristic 0.  While it would be very interesting to see analogues
of these results over other characteristics, the group theory involved
would be much more difficult.

\begin{thm}
  For all non-cyclic groups $G$, $$\frac{\mindeg \Q G}{\# G}\leq \frac
  32.$$
\end{thm}
\begin{proof}
  We split into 2 cases, depending on whether $G$ has a non-cyclic
  Sylow subgroup or not.
  
  \emph{Case 1.} Assume $S\subset G$ is a non-cyclic Sylow
  $q$-subgroup.  By Theorem \ref{P:pgroup}, $\mindeg \Q S\leq q(q+1)$.
  Since $|S|\geq q^2$,
  $$\frac{\mindeg \Q S}{\#S}\leq 1+\frac 1q\leq \frac 32.$$
  Furthermore, if $X_S$ and $Y_S$ are the sets we constructed in
  Proposition~\ref{P:pgroup} which
  realize this bound, $X_S$ has a fixed point and $Y_S$ does not, so
  we may apply Lemma \ref{T:ind} to see that $\mindeg \Q G\leq
  (G:S)q(q+1)$ so $$\mdoo\leq \frac{q(q+1)(G:S)}{\#G}\leq \frac 32.$$
  
  \emph{Case 2.} Now, assume that all Sylow subgroups of $G$ are
  cyclic. Such groups have been classified by H\"older, Burnside, and
  Zassenhaus \cite{dR96}.  They are exactly groups of the form
  \begin{equation}\label{eq:1}
    G=\inn{a,b:a^m=b^n=1, bab^{-1}=a^r}
  \end{equation}
  for some $m,n,r\in\Z$, where $r^n\equiv 1\pmod m$ and $m$ and
  $n(r-1)$ are coprime.
  
  Let $p$ be a prime dividing $m$, and let $G$ act on $\Fp$ by
  $^{a^{\ell}b^k}d=r^kd+\ell$. The number of fixed points of any
  element $a^\ell b^k\in G$ is simply the number of solutions to the
  equation
  \begin{equation*}
    (r^k-1)d=-\ell\pmod p.
  \end{equation*}
  Every element fixes an affine subspace of $\Fp$, that is, a set with
  0,1 or $p$ points.  If $g$ fixes $p$ points, then it is in the kernel
  $K$ of the action on $\Fp$.
  
  Consider the action of $G'=G/K$ on $\Fp$. Each nontrivial element of
  $G'$ fixes 0 or 1 elements of $\Fp$.
  Note that $b\in G$ fixes $0\in\Fp$, but no other element since $p$
  and $r-1$ are coprime.  Therefore, $G'$ is a non-regular Frobenius
  group (by our second definition).  Thus
  \begin{equation*}
    \mdoo\leq \frac{\mindeg\Q{G'}}{\#G'}\leq\frac 43.\qedhere
  \end{equation*}
\end{proof}
\begin{remark}
  As we mentioned before, we cannot hope for such a strong bound if our field
  $k$ has positive characteristic.  The first reduction step will
  imply that the bound is true for all groups except those with a
  non-cyclic Sylow $p$-subgroup, and all other Sylow subgroups cyclic,
  a much more complicated class of groups than those with all Sylow
  subgroups cyclic.
\end{remark}

Since $\mindeg\Q{\mathbb{F}_2\times \mathbb{F}_2}=6$, this bound is
sharp. However for most groups, it is actually quite bad.

For example, consider $A_4$, the alternating group of degree 4. As
before, we let
$K_4\lhd A_4$ denote the subgroup $K_4=\langle (12)(34),(13)(24)\rangle$.
One can check that the actions
\begin{align*}
X&\cong\frac{A_4}{A_4}\bigsqcup\frac{A_4}{\inn{(12)(34)}}\\
Y&\cong\frac{A_4}{K}\bigsqcup\frac{A_4}{\inn{(123)}}
\end{align*}
are \lineq over $\Q$ and of degree 7. A computer search shows that, in fact,
$\mindeg \Q{A_4}=7$, when our theorem implies it is $\leq 18$.

Here are a few results that give better bounds for certain classes of groups.

\begin{thm}\label{non-solv}
If $G$ is non-cyclic and $$\mdoo>\frac 34,$$ $G$ is solvable.
\end{thm}
\begin{proof}
Every non-abelian simple group has a 2-Sylow of order at least $8$
or subgroup isomorphic to $A_4$.  Thus $\mdoo\leq \frac 34$.  If
$G$ is not solvable, it has a non-abelian simple subquotient, so
$\mdoo\leq \frac 34$.
\end{proof}

\begin{thm}
$$\mdoo>\frac 43$$ if and only if $G\cong \mathbb{F}_2\times
\mathbb{F}_2$.
\end{thm}
\begin{proof}
The ``if'' direction is known.  For the ``only if,'' let $G$ be
such a group. If $G$ has only cyclic Sylow subgroups, or a
non-cyclic Sylow $q$-subgroup for $q>2$, then $\mdoo \leq \frac 43$, and if
the Sylow 2-subgroup of $G$ has order $>4$, $\mdoo\leq \frac 34$.  Thus,
$G$ must have a Sylow 2-subgroup $S\cong \mathbb{F}_2\times \mathbb{F}_2$.

By Theorem~\ref{non-solv}, $G$ is solvable.  Let $A$ be a minimal
normal abelian subgroup (which is thus elementary abelian).  If $A$ is
a $p$-group for $p>2$, then $G/A$ is not cyclic, and $\mdoo \leq
\frac{3}{2\#A}< \frac{4}{3}.$

Thus, $A$ is a $2$-group, and the quotient by it is cyclic.  Since
$G/A$ has a unique normal Sylow 2-subgroup, $G$ does as well.
By Schur-Zassenhaus, $G$ is a semi-direct product of $A$ and a
subgroup $A'$ which acts faithfully on $A$ (since the kernel of such
an action would be an abelian normal subgroup which is not a 2-group).

Since $\Aut(\mathbb{F}_2\times \mathbb{F}_2)\cong\mathbb{F}_3$, either
$\#A'=3$ and $G\cong A_4$, a case we have already ruled out, or $\#A=1$.
\end{proof}

For ``most'' groups, $\mdoo$ is smaller still.  While there are no
consistent results along these lines, we have found a number of first steps.

\begin{thm}
  If $\mdoo\geq 1$, then $G$ is of the form $C\ltimes
  \mathbb{F}_q^\epsilon$ where $C$ is cyclic, $\epsilon\in\{1,2\}$, and
  $q$ is a prime number.  Furthermore, $C$ acts faithfully on
  $\mathbb{F}_q^\epsilon$.
\end{thm}
\begin{proof}
Note that if $G$ has a non-cyclic quotient, then $\mdoo \leq 3/4$.
Thus all quotients of $G$ must be cyclic.  Similarly, by
Theorem~\ref{non-solv}, $G$ is solvable.

First, assume $G$ has non-cyclic Sylow $\ell$- and $q$-subgroups for
distinct primes $\ell\neq q$.  We may assume $G$ is solvable, so let
$A\lhd G$ be a normal abelian subgroup, and let $S$ be a Sylow
subgroup of $A$.  The subgroup $S$ is characteristic in $A$, and hence
normal in $G$, and $G/S$ must have a non-cyclic $\ell$- or $q$-Sylow.
Thus, we may exclude this case.

Now assume $G$ has a non-cyclic Sylow $q$-subgroup for exactly one prime $q$.
Let $A$ and $S$ be as above.  If $(\# S,q)=1$, then $G/S$ has a
non-cyclic Sylow subgroup.  Thus $S$ is a Sylow $q$-subgroup of $A$. Thus
$G/S$ is cyclic. In
particular, $G$ has a unique, normal Sylow $q$-subgroup$R$. Applying
Schur-Zassenhaus again, $G\cong C\ltimes R$, where $C\cong G/R$.

If $G$ has only cyclic Sylow subgroups, then $G$ is of the form
described in equation \eqref{eq:1}.  As we have already seen, if $m$
is not prime, this group has a non-cyclic quotient.  Thus
$\inn{a}\cong \mathbb{F}_q$, and obviously $G\cong \inn{b}\ltimes \inn{a}$.

In either case,  the map $C\to\Aut R$ must be
injective, since if $C'\subset C_G(R)$, $C'\lhd G$ and the
quotient $G/C'$ is not cyclic.  
\end{proof}

\begin{thm}
  For all $\epsilon>0$, there are only finitely many groups such that
  $\mdoo\geq 1+\epsilon$, but infinitely many such that $\mdoo>1$.
\end{thm}
\begin{proof}
  We need only prove the theorem for $\epsilon=1/n$ for $n\in \Z$.  By
  the characterization above, if $\mdoo\geq \frac{n+1}{n}$, then $G$
  corresponds to an element of $\Aut(\mathbb{F}_q^\epsilon)$ for some
  prime $q\leq n$ and $\epsilon\in\{1,2\}$.  These groups are finite,
  so there are only finitely many choices for $G$.

  We note that when $n$ is prime,
  $$\frac{\mindeg\Q{D_n}}{\#D_n}=\frac{n+1}{n}$$
  since in this case all
  proper subgroups of $D_n$ are cyclic, so the pair constructed in
  Section~\ref{sec:frobenius-groups} is the unique ``irreducible'' pair
  of non-isomorphic \lineq $G$-sets, in the sense that all others must
  this pair as a subset.  This exhibits infinitely many groups such that
  $\mdoo >1$.
\end{proof}

\section{Beaulieu's construction}\label{sec:beaul-constr}

In this section, we will describe the most fruitful known construction
of transitive $G$-sets, and how it allows to transfer some of our
results obtained in a highly non-transitive context to the transitive
case.

We let $X$ and $Y$ be $G$-sets, \lineq over $k$. By fixing bijections
$X\to\{1,\cdots,n\}$ and $Y\to\{1,\cdots,n\}$, where $n=|X|=|Y|$,
we obtain natural homomorphisms $\vp_X:G\to S_n$ and $\vp_Y:G\to S_n$,
where as usual, $S_n=\Sym(\{1,\ldots,n\})$.  The
homomorphisms $\vp_X$ and $\vp_Y$ are not unique, but any two
choices of $\vp_X$ will differ by an inner automorphism of $S_n$.
Fix bijections $\sigma_X:X\to \{1,\ldots,n\}$ and $\sigma_Y:Y\to \{1,\ldots,n\}$.

Let $X'$ be the $S_n$-set $S_n/\vp_X(G)$ (and similarly for $Y'$).
Note that isomorphism class of $X'$ is not changed by replacing
$\vp_X(G)$ with a conjugate, and thus will not depend on the
choice of $\vp_X$.

For simplicity, we will assume that the $G$-sets $X$ and $Y$ are
faithful (i.e., the homomorphisms $\vp_X$ and $\vp_Y$ are
injective).

\begin{thm}\label{T:beaump} If $G$-sets $X$ and $Y$ of degree
  $n$ are \lineq over $k$, then the $S_n$-sets $X'$ and $Y'$ are \lineq
  over $k$ as well.
\end{thm}
This theorem was originally proved for characteristic 0 by P. Beaulieu in her
Ph.D. thesis \cite{B91}.  When $p=0$, our proof essentially reduces
to a restatement of hers.
\begin{proof}
We will apply Lemma~\ref{L:restr} to $S_n$.  
Let $C$ be a \cmp subgroup of $S_n$.  If $C$ is not conjugate to a
subgroup of $\vp_{\res GC X}(G)=H_X$ or $\vp_{\res GC Y}(G)=H_Y$, then
$\Fix_{X'}(C)=\Fix_{Y'}(C)=\emptyset$

Thus, we need only consider \cmp subgroups contained in $H_X$ or $H_Y$.
If $C\subset H_X$, then there is a subgroup $K\subset G$ such that
$C=\vp_X(K)$.  Let $C'=\vp_Y(K)$. Now, consider the actions $\res GK X$
and $\res GK Y$. Since $X\lnq_k Y$, and $K$ is \cmp, there is an
isomorphism of $K$-sets $\tau:\res GK X\to\res GK Y$ by Lemma
\ref{L:restr}.

Thus, the permutation $s_K=\sigma_Y\circ \tau\circ \sigma_X^{-1}$ 
intertwines the action of $C$ and $C'$ on $\{1,\ldots,n\}$, that is, 
$s_KCs_K^{-1}=C'$.  

Now, consider the set $L_{C}(H_X)$.  This set can be partitioned as
\begin{equation*}
L_{C}(H_X)=\bigsqcup L_C^{C_1}(H_X)
\end{equation*}
as $C_1$ ranges over conjugates of $C$ contained
in $H$, and 
\begin{equation*}
  L_C^{C_1}(H_X)=\{g\in L_C(H_X)|g^{-1}Cg=C_1\}.
\end{equation*}
Now, note that $L_C^{C_1}(H_X)s_{K_1}=L_{C}^{C_1'}(H_Y)$,
where $\vp_X(K_1)=C_1,\vp_Y(K_1)=C_1'$.  This gives a bijection between
$L_C(H_X)$ and $L_C(H_Y)$. 

By Lemma \ref{L:fixcount}, this implies that
$\#\Fix_{X'}(C)=\#\Fix_{Y'}(C)$, and by Lemma \ref{L:restr}, we see that
$X'\lnq_k Y'$.
\end{proof}

The above theorem contains no information about whether $X'$ and $Y'$
are isomorphic as $S_n$-sets. This, of course, occurs exactly when
$H_X$ and $H_Y$ are conjugate in $S_n$. 

This is true if and only if the actions
$X$ and $Y$ are similar, that is, when there is a map $\mu:X\to Y$ and
an automorphism $\psi$ of $G$ such that
$\mu({^gx})={^{\psi(g)}\mu(x)}$.  This is a weaker condition than
requiring $X\cong_GY$, though these conditions are closely related.

This shows that the converse of Theorem \ref{T:beaump} is obviously
false, since we could take the actions $X$ and $Y$ to be similar but not
linearly equivalent over $k$.  In this case, $\vp_X(G)\sim\vp_Y(G)$ so
$X'$ and $Y'$ are isomorphic and thus \lineq over $k$.

\begin{question}
Are $X'$ and
$Y'$ are \lineq if and only if $X$ is \lineq to a $G$-set similar to $Y$?
\end{question}

This seems unlikely, since examples have been constructed by Perlis
\cite{rP77} of pairs of permutation groups $(G,X)$ and $(G',Y)$ such that $G$
and $G'$ are not isomorphic but $X'$ and $Y'$ are linearly equivalent
for a field of characteristic 0.

Since at present, we cannot answer the question above in general, let
us address a weaker form: When we can be sure that $k[X']\ncong
k[Y']$?

\begin{lem}\label{L:ncong}
  If $k[X']\cong k[Y']$ then for any \cmp subgroup $H_1\subset G$,
  there exists a subgroup $H_2\subset G$ such that $\res G{H_1}X$ and
  $\res G{H_2}Y$ are similar. In particular, 
  \begin{equation*}
    \#\Fix_X(H_1)=\#\Fix_Y(H_2).
  \end{equation*}
\end{lem}
While somewhat crude, this simple criterion allows us to construct
many pairs of $G$-sets for which we can be sure that $k[X']\ncong
k[Y']$.

\begin{proof}
Assume that $X'$ and $Y'$ are \lineq over $k$.  If $H_1$ is a \cmp
subgroup of $G$, then $\vp_X(H_1)$ must be conjugate to some \cmp
subgroup of $\vp_Y(G)$ by Lemma \ref{L:restr}, which must be of
the form $\vp_Y(H_2)$, for some $H_2\subset G$.  Translating back
into $G$-sets, this means that $\res G{H_1}X$ and $\res G{H_2}Y$
are similar.

Since any similarity of $G$-sets must preserve fixed points,
\begin{equation*}
\#\Fix_X(H_1)=\#\Fix_Y(H_2).\qedhere
\end{equation*}
\end{proof}

\section{An Application}\label{sec:an-application}

In general, it is quite difficult to find transitive $G$-sets
which are linearly equivalent but not isomorphic.  A number of
examples for small degrees have been studied by Perlis \cite{rP77},
Feit \cite{wF80}, Guralnick and Wales \cite{rG83,rG85}, and
DeSmit and Lenstra have recently given a more general construction
for $G$ solvable \cite{bD00}, but for the most part, this field
remains wide open.

Beaulieu's construction gives us a method of constructing a wide
variety of $G$-sets.  For example, it implies that the stabilizers of
\lineq $G$-sets have no special properties other than not being cyclic
mod $p$:
\begin{thm}\label{T:all-stabilizers}
  If $G$ is a group which is not \cmp, then for some $n$ there exist
  subgroups $G_1,G_2\subset S_n$ such that $G_1\cong G_2\cong G$ and
  $S_n/G_1$ and $S_n/G_2$ are \lineq over $k$ but not isomorphic.
\end{thm}
\begin{proof}
  The $G$-sets $\XM$ and $\YM$ defined in \eqref{eq:defXH} and
  \eqref{eq:defYH} where $\mc M$ contains the set of all maximal
  subgroups of $G$ and satisfies $\bigcap_{M_i\in \mc M}M_i=\{1\}$ are
  linearly equivalent by Theorem~\ref{T:exist}. Now,
  $G_1=\vp_{X_G}(G)$ and $G_2=\vp_{Y_G}(G)$ are precisely the
  subgroups we were looking for.  Since the $G$-sets are faithful, the
  maps are injective and thus isomorphisms onto their respective
  images.  By Theorem \ref{T:beaump}, $S_n/G_1\lnq_k S_n/G_2$.
  Lemma~\ref{L:ncong} shows that these are not isomorphic.
\end{proof}

Beaulieu's construction also allows us to show that $\mindeg k{S_n}$
is much smaller than the bounds shown in
Section~\ref{sec:bounding-degrees}, and that there is no lower bound
over all groups of $\frac{\mindeg k{G}}{\#G}$.  In fact, the examples
constructed by de Smit and Lenstra in \cite{bD00} show that no such
lower bound can be applied to the class of solvable groups, nilpotent
groups or
$q$-groups for any $q$.

\begin{thm}
For any field $k$, and any $\epsilon >0$, there exists an $N$ such that for all
$n>N$, $$\frac{\mindeg k{S_n}}{\#S_n}<\epsilon.$$

That is $$\lim_{n\to \infty}\frac{\mindeg k{S_n}}{\#S_n}=0.$$
\end{thm}
\begin{proof}
Fix a group $G$ which is not cyclic mod $q$ for any $q$.  Using
Theorem~\ref{T:all-stabilizers}, for some $N$, there are subgroups 
$H,H'\subset S_N$, which we can think of as subgroups of $S_n$ for any
$n\geq N$ by the standard inclusion maps, such that $H\cong H'\cong G$,
$S_n/H\ncong_{S_n}S_n/H'$, and $S_n/H\lnq_k S_n/H'$.  Thus
$$\frac{\mindeg k{S_n}}{\#S_n}\leq\frac{1}{\#G}$$
for all $n\geq N$.
Since there exist groups not cyclic mod any prime $q$ of arbitrary order (for
example $S_m$ as $m\geq 4$), the limit is
proved.
\end{proof}

In general, Beaulieu's construction helps us to turn non-transitive
constructions into transitive ones.  For example, I am not aware of
any previous example of a pair of transitive actions which are \lineq
over all fields but not isomorphic.  But since we can easily construct
non-transitive examples, Beaulieu's construction will now allow us to
construct as many of these as we would like.

For example, we can use the smallest group not cyclic mod any prime,
$G=D_6=\inn{a,b:a^6=b^2=(ab)^2}$. Consider the following $G$-sets
\begin{align}\label{eq:example}
X_{D_6}&=\frac G{\inn{a^2}}\bigsqcup\frac G{\inn{b}}\bigsqcup
  \frac G{\inn{ab}}\bigsqcup\frac G{\inn{a^3}}\bigsqcup2\cdot \left(\frac GG\right)\\
Y_{D_6}&=\frac G{\inn{a^2,b}}\bigsqcup\frac G{\inn{a^2,ab}}\bigsqcup
  \frac G{\inn{a}}\bigsqcup\frac G{\{1\}}\bigsqcup2\cdot \left(\frac G{\inn{b,a^3}}\right).
\end{align}

One can check that these $G$-sets are \lineq over any field, since
their restrictions to any proper subgroup are isomorphic.  Since
$|X_{D_6}|=|Y_{D_6}|=24$, we see that $S_{24}$ has transitive $G$-sets
linearly equivalent over any field.

In fact, this can be expanded further.  For any pair of $G$-sets,
there is a set of primes $\mc{P}_{X,Y}$, which is exactly the
primes $p$ such that $k[X]\ncong k[Y]$ for fields of
characteristic $p$.  This set is either all primes, or a finite
set dividing the order of $G$ (in fact, dividing the order of the
stabilizer of at least one point in $X$ or $Y$).

\begin{thm}\label{T:xdiv}
Given an arbitrary finite set of primes
$\mc{P}=\{p_1,\ldots,p_n\}$, there exists a group $G$, and a pair
of transitive $G$-sets $X$ and $Y$ such that $k[X]\cong k[Y]$ if
and only if $p=\char k\notin \mc{P}$, i.e. $\mc P= \mc P_{X,Y}$.
\end{thm}
\begin{proof}
We have already done the case where $\mc{P}=\{\}$.

Next we tackle singletons.  Let $\mc{P}=\{p\}$.  If $G$ is a
$p$-group, then for any $X$ and $Y$ \lineq over any field of
characteristic different from $p$, $\mc{P}_{X,Y}=\{p\}$ (since all
subgroups are \cmp).  By a construction of de Smit and Lenstra given
in \cite{bD00}, there exists a $p$-group $G_p$ with transitive
$G_p$-sets $X_p,Y_p$ \lineq over a field of characteristic 0 (and, in
fact, with degree $p^3$).

Now, for $\mc{P}$ general, the exterior Cartesian products
\begin{equation*}
X=\prod_{p\in\mc{P}}X_p\qquad\qquad Y=\prod_{p\in\mc{P}}Y_p
\end{equation*}
are obviously \lineq as $G=\prod_{p\in\mc{P}}G_p$-sets
for fields of characteristic $q\notin\mc{P}$.
If $p_i\in \mc{P}$, the $p$-subgroup
\begin{equation*}
H_{p_i}=\Stab_{G_{p_i}}(x_{p_i})\subset G
\end{equation*}
for $x_{p_i}\in X_{p_i}$ fixes points of $X$ but not of $Y$, and so by
Lemma \ref{L:restr}, $X$ and $Y$ are not linearly equivalent for any
field of characteristic $p_i$.  Thus $\mc{P}_{X,Y}=\mc{P}$.
\end{proof}

This proof only realized examples for groups which are nilpotent.  The
question of which sets appear for groups which not nilpotent (or more
generally, for groups which are indecomposable) appears more subtle,
but, in fact, the answer is the same.

A fairly limited number of sets (for the most part, singletons) have thus forth
come to light as $\mc P_{X,Y}$ for transitive $G$-sets $X$ and
$Y$, with $G$ indecomposable, but in fact, any finite set of primes can appear.

\begin{thm}\label{T:xdiv2}
Theorem \ref{T:xdiv} holds, with the additional assumption that $G\cong
S_n$ for some integer $n$.
\end{thm}
\begin{proof}
Let $G$ be as in the proof of Theorem \ref{T:xdiv}, let $\mc M$ be the
set of cyclic subgroups of $G$, and let $\XM$ and $\YM$ be as defined in
\eqref{eq:defXH} and \eqref{eq:defYH}.  These $G$ sets are \lineq over
$\Q$ and thus over all fields of characteristic $q\notin\mc P$, since
$q\nmid\#G$.  Thus, $\mc P_{\XM,\YM}\subset \mc P$.

Now, we apply Beaulieu's construction to $\XM$ and $\YM$, and denote the
corresponding $S_n$-sets by $X'$ and $Y'$ as before.   By
Theorem \ref{T:beaump}, $\mc P_{X',Y'}\subset\mc P$.

On the other hand, assume $p\in \mc P$, and let $K$ be any non-cyclic,
\cmp subgroup of $G$.  Let $\Fix_X (K)$ is the unique trivial orbit, while
$\Fix_Y (K)$ is empty.  Thus, by Lemma \ref{L:ncong}, $X'$ and $Y'$ are not
\lineq over any field of characteristic $p$.  So, $\mc P_{X',Y'}=\mc P$.
\end{proof}

The reader may wonder why we did not use the $G$-sets constructed
in the proof of Theorem~\ref{T:xdiv}.  In this case, it becomes
unclear whether one can apply Lemma~\ref{L:ncong}
at the end of the proof of Theorem \ref{T:xdiv2} to ensure that
$\mc{P}_{X',Y'}\subset\mc{P}_{X,Y}$. 

These constructions unfortunately tend to lead to $S_n$-sets of quite
enormous degrees.  For example, the example \eqref{eq:example} has
degree $24!/12$.  It would be very interesting to find other
constructions of $G$-sets isomorphic over a fixed set of fields, or
over other rings, which
would better live up to the title of this paper.

\providecommand{\bysame}{\leavevmode\hbox to3em{\hrulefill}\thinspace}
\providecommand{\MR}{\relax\ifhmode\unskip\space\fi MR }
\providecommand{\MRhref}[2]{%
  \href{http://www.ams.org/mathscinet-getitem?mr=#1}{#2}
}
\providecommand{\href}[2]{#2}

% \bibliography{thesis}

\begin{thebibliography}{NZM91}

\bibitem[Acc]{rA70}
Robert D.~M. Accola, \emph{Two theorems on {R}iemann surfaces with noncyclic
  automorphism groups.}, Proc. Amer. Math. Soc. \textbf{25} (1970), 598--602.
  \MR{41 \#3747}

\bibitem[BB]{BB04}
Werner Bley and Robert Boltje, \emph{Cohomological {M}ackey functors in number
  theory}, J. Number Theory \textbf{105} (2004), no.~1, 1--37. \MR{MR2032439
  (2004k:11169)}

\bibitem[Bea]{B91}
Patricia Beaulieu, \emph{{On a new construction of subgroups inducing
  isomorphic representations}}, Ph. D. thesis, Louisiana State University,
  1991.

\bibitem[Bol]{rB97}
Robert Boltje, \emph{Class group relations from {B}urnside ring idempotents},
  J. Number Theory \textbf{66} (1997), no.~2, 291--305. \MR{98i:20006}

\bibitem[CR]{cC87}
Charles~W. Curtis and Irving Reiner, \emph{Methods of representation theory.
  {V}ol. {I}{I}}, John Wiley \& Sons Inc., New York, 1987, With applications to
  finite groups and orders, A Wiley-Interscience Publication. \MR{88f:20002}

\bibitem[dS]{dS04}
Bart de~Smit, \emph{On arithmetically equivalent fields with distinct
  {$p$}-class numbers}, J. Algebra \textbf{272} (2004), no.~2, 417--424.
  \MR{MR2028064 (2005f:11252)}

\bibitem[dSL]{bD00}
B.~de~Smit and H.~W. Lenstra, Jr., \emph{Linearly equivalent actions of
  solvable groups}, J. Algebra \textbf{228} (2000), no.~1, 270--285.
  \MR{2001f:20069}

\bibitem[Fei]{wF80}
Walter Feit, \emph{Some consequences of the classification of finite simple
  groups}, The Santa Cruz Conference on Finite Groups (Univ. California, Santa
  Cruz, Calif., 1979), Amer. Math. Soc., Providence, R.I., 1980, pp.~175--181.
  \MR{82c:20019}

\bibitem[Glu]{dG81}
David Gluck, \emph{Idempotent formula for the {B}urnside algebra with
  applications to the $p$-subgroup simplicial complex}, Illinois J. Math.
  \textbf{25} (1981), no.~1, 63--67. \MR{82c:20005}

\bibitem[Gur]{rG83}
Robert~M. Guralnick, \emph{Subgroups inducing the same permutation
  representation}, J. Algebra \textbf{81} (1983), no.~2, 312--319.
  \MR{84j:20010}

\bibitem[GW]{rG85}
Robert~M. Guralnick and David~B. Wales, \emph{Subgroups inducing the same
  permutation representation. {I}{I}}, J. Algebra \textbf{96} (1985), no.~1,
  94--113. \MR{86m:20006}

\bibitem[NZM]{iN91}
Ivan Niven, Herbert~S. Zuckerman, and Hugh~L. Montgomery, \emph{An introduction
  to the theory of numbers}, fifth ed., John Wiley \& Sons Inc., New York,
  1991. \MR{91i:11001}

\bibitem[Per1]{rP77}
Robert Perlis, \emph{On the equation $\zeta \sb{{K}}(s)=\zeta \sb{{K}'}(s)$},
  J. Number Theory \textbf{9} (1977), no.~3, 342--360. \MR{56 \#5503}

\bibitem[Per2]{rP78}
\bysame, \emph{On the class numbers of arithmetically equivalent fields}, J.
  Number Theory \textbf{10} (1978), no.~4, 489--509. \MR{80c:12014}

\bibitem[Rob]{dR96}
Derek J.~S. Robinson, \emph{A course in the theory of groups}, second ed.,
  Springer-Verlag, New York, 1996. \MR{96f:20001}

\bibitem[Sol]{lS67}
Louis Solomon, \emph{The {B}urnside algebra of a finite group}, J.
  Combinatorial Theory \textbf{2} (1967), 603--615. \MR{35 \#5528}

\bibitem[ST]{hS00}
H.~M. Stark and A.~A. Terras, \emph{Zeta functions of finite graphs and
  coverings. {I}{I}}, Adv. Math. \textbf{154} (2000), no.~1, 132--195.
  \MR{2002f:11123}

\bibitem[Sun]{tS85}
Toshikazu Sunada, \emph{Riemannian coverings and isospectral manifolds}, Ann.
  of Math. (2) \textbf{121} (1985), no.~1, 169--186. \MR{86h:58141}

\end{thebibliography}
% \bibliographystyle{amsalpha}
\end{document}